\documentclass[draft]{article}

\usepackage[titletoc,title]{appendix}
\usepackage[cp1250]{inputenc}
\usepackage[IL2]{fontenc}
\usepackage{yfonts,fancyhdr}
\usepackage{a4wide}
\usepackage[english]{babel}
\usepackage{euscript}
\usepackage{amstext,amsbsy,amscd,amssymb}
\usepackage{amsmath,enumerate}
\usepackage{amsfonts}
\usepackage{mathrsfs,tensor,xy}
\usepackage{graphics}
\usepackage{microtype}
\include{diagxy}

\let\rarr=\rightarrow

\let\veps=\varepsilon
\let\mcal=\mathcal
\let\mfrak=\mathfrak
\let\eus=\EuScript

\def\Z{\mathbb{Z}}
\def\R{\mathbb{R}}
\def\C{\mathbb{C}}
\def\imag{{\rm i}}

\def\End{\mathop {\rm End} \nolimits}
\def\Hom{\mathop {\rm Hom} \nolimits}

\def\Ind{\mathop {\rm Ind} \nolimits}

\def\ad{\mathop {\rm ad} \nolimits}

\def\Diff{\mathop {\rm Diff} \nolimits}

\def\GL{\mathop {\rm GL} \nolimits}

\def\SL{\mathop {\rm SL} \nolimits}

\def\diag{\mathop {\rm diag} \nolimits}

\def\Sol{\mathop {\rm Sol} \nolimits}

\def\tr{\mathop {\rm tr} \nolimits}
\def\htt{\mathop {\rm ht} \nolimits}

\newcommand{\mP}{\mathbb{P}}

\newcommand{\mS}{\mathbb{S}}
\newcommand{\mV}{\mathbb{V}}

\newcommand{\mC}{\mathbb{C}}

\newcommand{\lC}{\mathcal{C}}

\newcommand{\mR}{\mathbb{R}}
\newcommand{\mW}{\mathbb{W}}

\newcommand{\gog}{\mathfrak{g}}

\newcommand{\gop}{\mathfrak{p}}

\newcommand{\gol}{\mathfrak{l}}

\newsymbol\squares 1003

\long\def\proof #1{\noindent \emph{Proof.}\ #1 \hfill $\squares$
\medskip}

\newcounter{num}[section]
\numberwithin{equation}{section}
\numberwithin{num}{section}

\long\def\definition #1 {\refstepcounter{num} \noindent {\bf
Definition \thenum.} #1

\medskip}

\long\def\theorem #1{\refstepcounter{num} \noindent {\bf Theorem
\thenum.} #1

\medskip}

\long\def\lemma #1{\refstepcounter{num}  \noindent {\bf Lemma
\thenum.} #1

\medskip}

\long\def\remark{\noindent {\bf Remark.}\ }

\newcommand*\riso{%
  \xrightarrow[]{\raisebox{-0.25em}{\smash{\ensuremath{\sim}}}}%
}

\makeatletter
\newcommand*\if@single[3]{%
  \setbox0\hbox{${\mathaccent"0362{#1}}^H$}%
  \setbox2\hbox{${\mathaccent"0362{\kern0pt#1}}^H$}%
  \ifdim\ht0=\ht2 #3\else #2\fi
  }
\newcommand*\rel@kern[1]{\kern#1\dimexpr\macc@kerna}
\newcommand*\widebar[1]{\@ifnextchar^{{\wide@bar{#1}{0}}}{\wide@bar{#1}{1}}}
\newcommand*\wide@bar[2]{\if@single{#1}{\wide@bar@{#1}{#2}{1}}{\wide@bar@{#1}{#2}{2}}}
\newcommand*\wide@bar@[3]{%
  \begingroup
  \def\mathaccent##1##2{%
    \if#32 \let\macc@nucleus\first@char \fi
    \setbox\z@\hbox{$\macc@style{\macc@nucleus}_{}$}%
    \setbox\tw@\hbox{$\macc@style{\macc@nucleus}{}_{}$}%
    \dimen@\wd\tw@
    \advance\dimen@-\wd\z@
    \divide\dimen@ 3
    \@tempdima\wd\tw@
    \advance\@tempdima-\scriptspace
    \divide\@tempdima 10
    \advance\dimen@-\@tempdima
    \ifdim\dimen@>\z@ \dimen@0pt\fi
    \rel@kern{0.6}\kern-\dimen@
    \if#31
      \overline{\rel@kern{-0.6}\kern\dimen@\macc@nucleus\rel@kern{0.4}\kern\dimen@}%
      \advance\dimen@0.4\dimexpr\macc@kerna
      \let\final@kern#2%
      \ifdim\dimen@<\z@ \let\final@kern1\fi
      \if\final@kern1 \kern-\dimen@\fi
    \else
      \overline{\rel@kern{-0.6}\kern\dimen@#1}%
    \fi
  }%
  \macc@depth\@ne
  \let\math@bgroup\@empty \let\math@egroup\macc@set@skewchar
  \mathsurround\z@ \frozen@everymath{\mathgroup\macc@group\relax}%
  \macc@set@skewchar\relax
  \let\mathaccentV\macc@nested@a
  \if#31
    \macc@nested@a\relax111{#1}%
  \else
    \def\gobble@till@marker##1\endmarker{}%
    \futurelet\first@char\gobble@till@marker#1\endmarker
    \ifcat\noexpand\first@char A\else
      \def\first@char{}%
    \fi
    \macc@nested@a\relax111{\first@char}%
  \fi
  \endgroup
}
\makeatother

\newcommand\rsmraise[1]{%
  \ifx#1\displaystyle .8\else
    \ifx#1\textstyle .8\else
      \ifx#1\scriptstyle .6\else
        .45%
      \fi
    \fi
  \fi}

\title{Projective structure, $\widetilde{\mathrm{SL}}(3,\mR)$ and the symplectic Dirac operator}

\author{Marie Holíková, Libor Křižka, Petr Somberg}

\AtEndDocument{\bigskip{\footnotesize%
  (M.\,Holíková) \textsc{Department of Mathematics and Mathematical Education, Faculty of Education,\\ Charles University,
  Magdal\'eny Rettigov\'e 4, 116\,39 Praha 1, Czech Republic}\par
  \textit{E-mail address}: \texttt{marie.holikova@pedf.cuni.cz}\par
  \addvspace{\medskipamount}
  (L.\,Křižka) \textsc{Mathematical Institute of Charles University, Sokolovsk\'a 83, 180\,00 Praha 8, Czech Republic} \par
  \textit{E-mail address}: \texttt{krizka@karlin.mff.cuni.cz} \par
  \addvspace{\medskipamount}
  (P.\,Somberg) \textsc{Mathematical Institute of Charles University, Sokolovsk\'a 83, 180\,00 Praha 8, Czech Republic} \par
  \textit{E-mail address}: \texttt{somberg@karlin.mff.cuni.cz}
}}

\begin{document}
\date{}
\maketitle

\begin{abstract}
Inspired by the results on symmetries of the symplectic Dirac operator,
we realize symplectic spinor fields and the symplectic Dirac operator in
the framework of (the double cover of) homogeneous projective structure
in two real dimensions.
The symmetry group of the homogeneous model of the double cover of projective
geometry in two real dimensions is $\smash{\widetilde{\SL}}(3,\mR)$.


\end{abstract}

\thispagestyle{empty}

\tableofcontents


\section*{Introduction}
\addcontentsline{toc}{section}{Introduction}

There are symplectic counterparts of the notions of the spinor field
and the Dirac operator on manifolds with metaplectic structure, see e.g.\ \cite{Kostant1974},
\cite{Habermann-Habermann2006}. The aim of the present short article
is to describe a realization of symplectic spinor fields and the symplectic
Dirac operator $D_s$ in the framework of (the double covering of) the
geometry of projective structure in the real dimension two. An inspiration
for this development comes from the recent results on symmetries of the
symplectic Dirac operator, \cite{DeBie-Holikova-Somberg2015}.

We shall briefly comment on the content of our article. In Section \ref{sec:symmetr} we
start with a motivation for our work, namely the question of differential
symmetries of the symplectic Dirac operator in real dimension two. Some
information on this structure already appeared in \cite{DeBie-Holikova-Somberg2015}, but its
meaning and interpretation was not clear at that time. Here we find a natural explanation in
terms of the Lie algebra $\mfrak{sl}(3,\mR)$ as an organizing principle for the
symmetry algebra of the solution space of $D_s$. To this aim, we introduce
in Section \ref{sec:projstr} the homogeneous projective structure in the real dimension
two and describe its basic geometrical and representation theoretical
properties. In Section \ref{sec:singular vectors} we briefly review the
procedure called F-method, which is applied in the last Section \ref{sec:sl3sdo}
to the simple metaplectic components of the Segal-Shale-Weil representation
(twisted by a character of the central generator of the Levi factor) as an inducing
representation for generalized Verma modules associated to the Lie algebra
$\mfrak{sl}(3,\mR)$ and its maximal parabolic subalgebra. In this way we produce the
symplectic Dirac operator as an ${\widetilde{\SL}}(3,\mR)$-equivariant differential operator
on the double covering of the real projective space $\mathbb{RP}^2$.

Let us highlight the meaning of the solution space for the
symplectic Dirac operator $D_s$, regarded as an equivariant differential
operator acting on the $\smash{\widetilde{\SL}}(3,\mR)$-principal series representation induced
from the character twisted metaplectic representation of $\smash{\widetilde{\SL}}(2,\mR)$.
The solution space of $D_s$ on the symplectic space $\mR^2$ was already determined in
\cite{DeBie-Somberg-Soucek2014}, as a consequence of
the metaplectic Howe duality for the pair $(\mathfrak{mp}(2,\mR), \mathfrak{sl}(2,\mR))$. The
underlying Harish-Chandra module of $\ker D_s$ for the Harish-Chandra pair
$(\mathfrak{sl}(3,\mR),\mathrm{SU}(2))$ with the maximal compact subgroup
$\mathrm{SU}(2)\subset\smash{\widetilde{\SL}}(3,\mR)$ is a unitarizable irreducible representation,
equivalent to the exceptional representation of $\smash{\widetilde{\SL}}(3,\mR)$
which is associated with the minimal coadjoint orbit, cf.\ \cite{Oersted2000}, \cite{Torasso1983}.

As for the notation used throughout the article, we consider the pair $(G,P)$ consisting
of a connected real reductive Lie group $G$ and its parabolic subgroup $P$. In the Levi
decomposition $P=LU$, $L$ denotes the Levi subgroup and $U$ the unipotent subgroup of $P$.
We write $\gog(\mR)$, $\gop(\mR)$, $\gol(\mR)$, $\mfrak{u}(\mR)$ for the
real Lie algebras and $\gog$, $\gop$, $\gol$, $\mfrak{u}$ for the complexified Lie
algebras of
$G$, $P$, $L$, $U$, respectively. The symbol $U$ applied to a Lie
algebra denotes its universal enveloping algebra, and similarly
$\widetilde{\phantom{nn}}$ applied
to a Lie group denotes its double cover.

\section{Symmetries of the symplectic Dirac operator}
\label{sec:symmetr}

In the present section we start with a motivation for our further
considerations. Let $(\mR^2,\omega)$ be the $2$-dimensional 
symplectic vector space with the canonical symplectic form $\omega=dx\wedge dy$,
where $(x,y)$ are the canonical linear coordinate functions on $\R^2$.
We denote by $\partial_x, \partial_y$ the coordinate vector
fields and by $e_1,e_2$ corresponding symplectic frame fields 
acting on a symplectic spinor
$\varphi\in \C[\R^2] \otimes_\C \mathcal{S}(\mR)$ by
\begin{align}
e_1 \cdot\varphi= \imag q \varphi,\qquad e_2 \cdot\varphi=\partial_q \varphi.
\end{align}
Here $\mathcal{S}(\mR)$ is the Schwartz space of rapidly decreasing complex functions on
$\mR$ equipped with the coordinate function $q$.

The basis elements $\{{X},{Y},{H}\}$ of the metaplectic
Lie algebra $\mathfrak{mp}(2,\mR)$, the double
cover of the symplectic Lie algebra $\mathfrak{sp}(2,\mR)\simeq \mathfrak{sl}(2,\mR)$, act on the
function space $\C[\R^2] \otimes_\C \mathcal{S}(\mR)$ by
\begin{align}\label{HXYxy}
 X=- y \partial_x-{\textstyle {\imag \over 2}} q^2,\qquad H=- x \partial_x+y \partial_y+q \partial_q+{\textstyle {1 \over 2}}, \qquad Y=- x \partial_y-{\textstyle {\imag \over 2}} \partial_q^2
\end{align}
and satisfy the commutation relations
\begin{align}
[H,X]=2X, \qquad [X,Y]=H, \qquad [H,Y]=-2Y.
\end{align}
The operators \eqref{HXYxy} preserve the homogeneity in the variables $x,y$,
and the $\mathfrak{mp}(2,\mR)$-equivariant operators
\begin{align}\label{sl2howedual}
 X_s = y \partial_q + \imag x q,\qquad  E  = x\partial_x+y\partial_y+{\textstyle {1 \over 2}}, \qquad
 D_s = \imag q \partial_y - \partial_x\partial_q
\end{align}
form a Lie algebra isomorphic to $\mathfrak{sl}(2,\mR)$. The operator $D_s$ is termed
the symplectic Dirac operator on $\R^2$.

A differential operator $P$ is called a symmetry of $D_s$
provided there exists another differential operator $P'$
such that $P'D_s=D_sP$. Consequently, the symmetry differential
operators preserve the solution space of $D_s$. The vector space of first
order (in all variables $x,y,q$) symmetries was described
in \cite{DeBie-Holikova-Somberg2015}, and there is a Lie algebra structure given by the
commutator of two symmetry differential operators. The commutators of
elements in \eqref{HXYxy} with the symmetry differential operators 
classified in \cite{DeBie-Holikova-Somberg2015} yield the following 
result, whose proof is a straightforward
but tedious computation.
\medskip

\lemma{The solution space of the symplectic Dirac operator on $\R^2$ is preserved by the following differential operators:
\begin{enumerate}
\item[1)]
The couple of commuting differential operators
\begin{align}
\mcal{O}_1&=x^2\partial_x+xy\partial_y-{\textstyle \frac{1}{2}}xq\partial_q+{\textstyle \frac{\imag}{2}}y\partial_q^2+{\textstyle \frac{1}{2}}x
=-{\textstyle \frac{1}{2}}xH-yY+{\textstyle \frac{1}{2}}x E+{\textstyle \frac{1}{2}}x, \\
\mcal{O}_2&=xy\partial_x+y^2\partial_y+{\textstyle \frac{1}{2}}yq\partial_q+{\textstyle \frac{\imag}{2}}xq^2+y
 ={\textstyle \frac{1}{2}}yH-xX+{\textstyle \frac{1}{2}}y E+{\textstyle \frac{1}{2}}y
\end{align}
satisfies
\begin{equation}
[D_s, \mcal{O}_1]={\textstyle \frac{3}{2}}xD_s
\qquad \text{and} \qquad
[D_s, \mcal{O}_2]={\textstyle \frac{3}{2}}yD_s, 
\end{equation}
which is equivalent to 
$D_s\mcal{O}_1=\big(\mcal{O}_1+\frac{3}{2}x\big)D_s$ and
$D_s \mcal{O}_2=\big(\mcal{O}_2+\frac{3}{2}y\big)D_s$, respectively.
The operators $\mcal{O}_1, \mcal{O}_2$ increase the homogeneity in the variables $x,y$ by one.

\item[2)]
The couple of commuting differential operators
\begin{equation}
\partial_x\quad \text{and} \quad \partial_y
\end{equation}
commutes with $D_s$ and decreases the homogeneity
in the variables $x,y$ by one.
\end{enumerate}}

The operators $\mcal{O}_1,\mcal{O}_2$ turn out to be useful in the
light of the following observation,
whose proof is again straightforward and left to the reader.
\medskip

\theorem{The first order symmetry differential operators
\begin{align}\label{symoper}
\{\partial_x,\partial_y,{H},{X},{Y},E,\mcal{O}_1,\mcal{O}_2\}
\end{align}
of the symplectic Dirac operator
$D_s$ fulfill the following non-trivial commutation relations
\begin{align}
\begin{array}{lll}
\, [\partial_x,\mcal{O}_2] =-X, & \quad & [\partial_y,\mcal{O}_1] =-Y, \\
\, [\partial_y,\mcal{O}_2] = \frac{1}{2}(3E+H), & &[\partial_x,\mcal{O}_1] = \frac{1}{2}(3E-H), \\
\, [\mcal{O}_2,H]=-\mcal{O}_2, & &[\mcal{O}_1,H]=\mcal{O}_1, \\
\, [\mcal{O}_2,E]=-\mcal{O}_2, & &[\mcal{O}_1,E]=-\mcal{O}_1, \\
\, [\mcal{O}_2,Y]=\mcal{O}_1,  & &[\mcal{O}_1,X]=\mcal{O}_2, \\
\, [\partial_x,{H}]=-\partial_x,  & & [\partial_y,{H}]=\partial_y,  \\
\, [\partial_x,{Y}]=-\partial_y,  & & [\partial_y,{X}]=-\partial_x, \\
\, [\partial_x,E]=\partial_x,  & & [\partial_y,E]=\partial_y, \\
\, [H,X]=2X, & & [H,Y]=-2Y, \\
\, [X,Y]=H.
\end{array}
\end{align}
The homomorphism of Lie algebras
\begin{align}
\langle \partial_x,\partial_y,
{H},{X},{Y},E,\mcal{O}_1,\mcal{O}_2 \rangle \rarr \mathfrak{sl}(3,\mR)
\end{align}
given by
\begin{align}
\begin{gathered}
-\partial_x\mapsto \begin{pmatrix}
0 & 0 & 0 \\
1 & 0 & 0 \\
0 & 0 & 0
\end{pmatrix}\!, \qquad
-\partial_y\mapsto \begin{pmatrix}
0 & 0 & 0 \\
0 & 0 & 0 \\
1 & 0 & 0
\end{pmatrix}\!,\\
 \mcal{O}_1\mapsto \begin{pmatrix}
0 & 1 & 0 \\
0 & 0 & 0 \\
0 & 0 & 0
\end{pmatrix}\!, \qquad
\mcal{O}_2\mapsto \begin{pmatrix}
0 & 0 & 1 \\
0 & 0 & 0 \\
0 & 0 & 0
\end{pmatrix}\!, \\
X \mapsto \begin{pmatrix}
0 & 0 & 0 \\
0 & 0 & 1 \\
0 & 0 & 0
\end{pmatrix}\!, \quad
H \mapsto \begin{pmatrix}
0 & 0 & 0 \\
0 & 1 & 0 \\
0 & 0 & -1
\end{pmatrix}\!, \quad
Y \mapsto \begin{pmatrix}
0 & 0 & 0 \\
0 & 0 & 0 \\
0 & 1 & 0
\end{pmatrix}\!,\\
E\mapsto \begin{pmatrix}
\frac{2}{3} & 0 & 0 \\
0 & -\frac{1}{3} & 0 \\
0 & 0 & -\frac{1}{3}
\end{pmatrix}
\end{gathered}
\end{align}
is an isomorphism of Lie algebras.}

As will be proved in \cite{Somberg-Silhan2016} by the techniques of 
tractor calculus, the operators \eqref{symoper} in the
variables $x,y,q$ are in fact all first order symmetry differential 
operators of $D_s$ in the base variables $x,y$.

In the remaining part of our article we interpret the results of
the present section in the framework of (the double
cover of) the projective structure in real dimension $2$.


\section{Generalized Verma modules and singular vectors}
\label{sec:singular vectors}

It is well-known that the $G$-equivariant differential operators acting on
principal series representations for $G$ can be recognized in the study
of homomorphisms between generalized Verma modules for the Lie algebra
$\mathfrak{g}$. The latter homomorphisms are determined by
the image of the highest weight vectors, referred to as the
singular vectors and characterized as being annihilated by the positive
nilradical $\mfrak{u}$.

An approach to find precise positions of singular vectors in the
representation space can be found in \cite{koss}, \cite{Krizka-Somberg2015},
\cite{Krizka-Somberg2016}.
Let $\mV$ denote a complex simple
highest weight $L$-module, extended to $P$-module by $U$ acting trivially.
We denote by $\mV^*$ the (restricted) dual $P$-module to $\mV$. Any
character $\lambda \in \Hom_P(\mfrak{p},\C)$ yields a $1$-dimensional 
representation $\C_\lambda$ of $\mfrak{p}$ by
\begin{align}
  Xv=\lambda(X)v,\quad X \in \mfrak{p},\, v \in \C .
\end{align}
Moreover, assuming that $\lambda \in \Hom_P(\mfrak{p},\C)$ 
defines a group character $e^\lambda \colon P \rarr \GL(1,\C)$ of $P$ 
and denoting by $\rho \in \Hom_P(\mfrak{p},\C)$ the character
\begin{align}
  \rho(X) = {\textstyle {1 \over 2}} \tr_\mfrak{u}\ad(X),
	\quad X \in \mfrak{p}, \label{eq:rho vector}
\end{align}
we may introduce a twisted $P$-module
$\mathbb{V}_{\lambda+\rho}\simeq \mathbb{V}\otimes_\C \mC_{\lambda+\rho}$ (with a twist
$\lambda +\rho$) where $p \in P$ acts on $v \in \mathbb{V}_{\lambda+\rho} \simeq \mathbb{V}$ 
(the isomorphism of vector spaces) by $e^{\lambda+\rho}(p)p.v$. In the rest of our article, 
$\mV=\mathbb{S}$ is one of the simple metaplectic
submodules of the Segal-Shale-Weil representation twisted by character
of $\GL(1,\R)_+$.

By abuse of notation, we call the highest weight modules induced from infinite-dimensional
simple $L$-modules generalized Verma modules, though they are not the objects
of the parabolic BGG category ${\fam2 O}^\gop$ because of the lack of $L$-finiteness
condition.
However, most of structural results required in the present article carry over to
this class of modules, cf.\ \cite{Krizka-Somberg2016}.

In general, for a chosen principal series representation of ${G}$ on the vector space
$\Ind_{{P}}^{{G}}(\mV_{\lambda+\rho})$ of smooth sections of the
homogeneous vector bundle
${G} \times_{{P}} \mV_{\lambda+\rho} \rarr {G}/{P}$
associated to a ${P}$-module $\mV_{\lambda+\rho}$, we compute the infinitesimal action
\begin{align}
\pi_\lambda \colon \mfrak{g} \rarr \mcal{D}(U_e) \otimes_\C \End \mathbb{V}_{\lambda+\rho}.
\end{align}
Here $\mcal{D}(U_e)$ denotes the $\C$-algebra of smooth complex linear
differential operators on
$U_e=\widebar{U}{P} \subset {G}/{P}$ ($\widebar{U}$
is the Lie group whose Lie algebra is the opposite nilradical $\widebar{\mathfrak{u}}(\R)$
to $\mathfrak{u}(\R)$), on the vector space
$\mcal{C}^\infty(U_e) \otimes_\C\! \mathbb{V}_{\lambda+\rho}$ of
$\mathbb{V}_{\lambda+\rho}$-valued smooth functions on $U_e$ in the non-compact
picture of the induced representation.

The dual vector space $\mcal{D}'_o(U_e)\otimes_\C\! \mV_{\lambda+\rho}$ of
$\mathbb{V}_{\lambda+\rho}$-valued distributions on $U_e$ supported on the
unit coset $o=e {P}\in {G}/{P}$ is
$\mcal{D}(U_e) \otimes_\C \End \mathbb{V}_{\lambda+\rho}$-module,
and there is an $U(\gog)$-module isomorphism
\begin{align}\label{vermadistr}
\Phi_\lambda \colon M^\gog_\gop(\mV_{\lambda-\rho})\equiv U(\gog)\otimes_{U(\gop)}\!\mV_{\lambda-\rho} \rarr  \mcal{D}'_o(U_e)\otimes_\C\! \mV_{\lambda+\rho} \simeq \eus{A}^\mfrak{g}_{\widebar{\mfrak{u}}}/I_e \otimes_\C\! \mathbb{V}_{\lambda+\rho}.
\end{align}
The exponential map allows to identify $U_e$ with the nilpotent Lie algebra
$\widebar{\mfrak{u}}(\R)$. If we denote by $\eus{A}^\mfrak{g}_{\widebar{\mfrak{u}}}$
the Weyl algebra of the complex vector space $\widebar{\mfrak{u}}$, then the vector
space $\mcal{D}'_o(U_e)$ can be identified as an
$\eus{A}^\mfrak{g}_{\widebar{\mfrak{u}}}$-module with the quotient of
$\eus{A}^\mfrak{g}_{\widebar{\mfrak{u}}}$ by the left ideal $I_e$ generated
by all polynomials on $\widebar{\mathfrak{u}}$ vanishing at the origin.

Let $(x_1,x_2,\dots,x_n)$ be the linear coordinate functions on $\widebar{\mfrak{u}}$ and
$(y_1,y_2,\dots,y_n)$ be the dual linear coordinate functions on $\widebar{\mfrak{u}}^*$.
Then the algebraic Fourier transform
\begin{align}
  \mcal{F} \colon \eus{A}^\mfrak{g}_{\widebar{\mfrak{u}}} \rarr \eus{A}^\mfrak{g}_{\widebar{\mfrak{u}}^*}
\end{align}
is given by
\begin{align}
  \mcal{F}(x_i) = -\partial_{y_i}, \qquad \mcal{F}(\partial_{x_i}) = y_i
\end{align}
for $i=1,2,\dots,n$, and leads to a vector space isomorphism
\begin{align}\label{eqn:FT}
\begin{gathered}
\tau \colon  \eus{A}^\mfrak{g}_{\widebar{\mfrak{u}}}/I_e \simeq \C[\widebar{\mfrak{u}}^*]  \riso  \eus{A}^\mfrak{g}_{\widebar{\mfrak{u}}^*}/\mcal{F}(I_e) \simeq \C[\widebar{\mfrak{u}}^*], \\
  Q\ {\rm mod}\ I_e  \mapsto \mcal{F}(Q)\ {\rm mod}\ \mcal{F}(I_e)
\end{gathered}
\end{align}
for $Q \in \eus{A}^\mfrak{g}_{\widebar{\mfrak{u}}}$. The composition of
\eqref{vermadistr} and \eqref{eqn:FT} gives the vector space isomorphism
\begin{align}\label{eqn:phiPM}
 \tau \circ \Phi_\lambda \colon U(\mfrak{g}) \otimes_{U(\mfrak{p})} \! \mathbb{V}_{\lambda-\rho} \riso \mcal{D}'_o(U_e) \otimes_\C \!\mV_{\lambda+\rho} \riso
\C[\widebar{\mathfrak{u}}^*]\otimes_\C\! \mV_{\lambda+\rho},
\end{align}
thereby inducing the $\mfrak{g}$-module action $\hat{\pi}_\lambda$ on
$\C[\widebar{\mathfrak{u}}^*]\otimes_\C\! \mV_{\lambda+\rho}$.
\medskip

\definition{Let $\mV$ be a complex simple highest weight $L$-module,
extended to a $P$-module by $U$ acting trivially. We define the $L$-module
\begin{align}
M_\gop^\gog(\mV)^{\mathfrak{u}}= \{v\in M^\gog_\gop(\mV);\, Xv=0\
\text{for all}\ X\in \mathfrak{u}\},
\end{align}
which we call the vector space of singular vectors.}

The vector space of singular vectors is for an infinite-dimensional complex simple
highest weight $P$-module $\mV$ an infinite-dimensional $L$-module. In the case
when $M_\gop^\gog(\mV)^\mathfrak{u}$ is a completely reducible $L$-module,
we denote by $\mW$ one of its simple $L$-submodule. Then we obtain the
$U(\gog)$-module homomorphism from $M_\gop^\gog(\mW)$ to $M_\gop^\gog(\mV)$, and moreover we have
\begin{align}
\Hom_{(\gog,P)}(M_\gop^\gog(\mW), M_\gop^\gog(\mV))
\simeq
\Hom_L(\mW, M_\gop^\gog(\mV)^\mathfrak{u}).
\end{align}
We introduce the $L$-module
\begin{align} \label{eqn:sol2l}
\Sol(\mathfrak{g},\mathfrak{p};\C[\widebar{\mathfrak{u}}^*]\otimes_\C\! \mV_{\lambda+\rho})^\mcal{F} =\{f \in \C[\widebar{\mathfrak{u}}^*]\otimes_\C\! \mV_{\lambda+\rho};\, \hat{\pi}_\lambda(X) f = 0\ \text{for all}\ X\in\mathfrak{u}\},
\end{align}
and by \eqref{eqn:phiPM}, there is an $L$-equivariant isomorphism
\begin{equation}
\label{eqn:phi}
\tau \circ \Phi_\lambda \colon M_\mathfrak{p}^\mathfrak{g}(\mV_{\lambda-\rho})^\mathfrak{u} \riso \Sol(\mathfrak{g},\mathfrak{p}; \C[\widebar{\mathfrak{u}}^*]\otimes_\C\! \mV_{\lambda+\rho})^\mcal{F}.
\end{equation}
The action of $\hat{\pi}_\lambda(X)$ on $\C[\widebar{\mathfrak{u}}^*]\otimes_\C\! \mV_{\lambda+\rho}$
produces a system of partial differential equations for the elements in
$\Sol(\mathfrak{g},\mathfrak{p}; \C[\widebar{\mathfrak{u}}^*]\otimes_\C\! \mV_{\lambda+\rho})^\mcal{F}$, which makes it possible to describe its structure completely in particular cases of interest as the
solution space of the systems of partial differential equations.

The formulation above has the following classical dual statement
(cf.\ \cite{Collingwood-Shelton1990} for the standard formulation in the category of
finite dimensional inducing $P$-modules, or \cite{Krizka-Somberg2016} for its
extension to inducing modules with infinitesimal character), which
explains the relationship between the geometrical problem of finding
$G$-equivariant differential operators between induced representations
and the algebraic problem of finding homomorphisms between generalized
Verma modules. Let $\mV$ and $\mW$ be two simple highest weight $P$-modules.
Then the vector space of $G$-equivariant differential operators
$\Hom_{\Diff(G)}(\Ind_P^G(\mV),\Ind_P^G(\mW))$
is isomorphic to the vector space of $(\gog,P)$-homomorphisms
$\Hom_{(\gog,P)}(M^\gog_\gop(\mW^*),M^\gog_\gop(\mV^*))$.


\section{The homogeneous projective structure in dimension $2$}
\label{sec:projstr}

A projective structure on a smooth manifold $M$ of
real dimension $n\geq 2$ is a class $[\nabla]$ of projectively
equivalent torsion-free connections, which define the same family of unparametrized
geodesics. A connection is projectively flat if and only if it is locally
equivalent to a flat connection. Given any non-vanishing
volume form $\omega$ on $M$, there exists a unique connection
in the projective class such that $\nabla\omega=0$. In the case
$n>2$ ($n=2$), the vanishing of the Weyl curvature tensor (the
Cotton curvature tensor) for $\nabla$ is equivalent to the projective
flatness of $[\nabla]$ and consequently, the existence of a local isomorphism
with the flat model of $n$-dimensional projective geometry on $\mathbb{RP}^n$
equipped with the flat projective structure given by the absolute
parallelism.

The homogeneous (flat) model of projective geometry in the real
dimension $2$ is $\mathbb{RP}^2\simeq G/P$, where $G$ is the connected 
simple real Lie group $\SL(3,\mR)$ and $P\subset G$ the parabolic
subgroup stabilizing the line $[v]\in\mR^3$ generated by a non-zero vector $v$ in
the defining representation $\mR^3$ of $G$.
Although our construction of equivariant differential operators is local,
the passage to the generalized flag manifold and sections of associated vector
bundles induced from half integral modules (e.g.,\ the simple metaplectic submodules
of the Segal-Shale-Weil representation) on $G/P$ requires the double (universal) cover
$\smash{\widetilde{G}}=\smash{\widetilde{\SL}}(3,\R)$ and its parabolic subgroup 
$\smash{\widetilde{P}}$. The Lie
group $\smash{\widetilde{\SL}}(3,\R)$ acts transitively on $S^2\simeq \mC\mP^1$,
the double (universal) cover of $\mathbb{RP}^2$, with parabolic stabilizer
$\smash{\widetilde{P}}=(\GL(1,\R)_+\times \smash{\widetilde{\SL}}(2,\R))\ltimes\mR^2$.
We notice that the double (universal) cover 
$\smash{\widetilde{\SL}}(3,\R)/\smash{\widetilde{P}}\simeq S^2\simeq \mC\mP^1$ 
is a symplectic manifold,
while $\mR\mP^2$ is non-orientable and hence not symplectic.

The questions discussed in our article can be treated by infinitesimal methods,
and so we introduce the complexified
Lie algebra $\mathfrak{g}=\mathfrak{sl}(3,\C)$ of $G$ and the Cartan
subalgebra $\mathfrak{h}\subset\mathfrak{g}$ by
\begin{align}
  \mathfrak{h}=\{\diag(a_1, a_2, a_3);\, a_1+a_2+a_3=0,\, a_1,a_2,a_3\in\C\}.
\end{align}
For $i=1,2,3$, we define $\varepsilon_i \in \mathfrak{h}^*$ by
$\varepsilon_i(\diag(a_1, a_2, a_3))=a_i$. The root system
of $\mathfrak{g}$ with respect to $\mathfrak{h}$ is
$\Delta = \{\pm (\veps_i - \veps_j);\, 1 \leq i < j \leq 3\}$,
the positive root system is $\Delta^+=\{\veps_i-\veps_j;\, 1\leq i<j \leq 3\}$
with the subset of simple roots $\Pi=\{\alpha_1,\alpha_2\}$,
$\alpha_1=\veps_1-\veps_2$, $\alpha_2=\veps_2-\veps_3$, and the fundamental weights
are $\omega_1=\veps_1$, $\omega_2=\veps_1+\veps_2$.
The subset $\Sigma=\{\alpha_2\}$ of $\Pi$ generates a root
subsystem $\Delta_\Sigma\subset\mfrak{h}^*$, and we associate to $\Sigma$ the
standard parabolic subalgebra $\mfrak{p}$ of $\mfrak{g}$ with
$\mfrak{p} = \mfrak{l} \oplus \mfrak{u}$.
The reductive Levi subalgebra $\mfrak{l}$ of $\mfrak{p}$ is
\begin{align}
\mfrak{l}= \mfrak{h} \oplus \bigoplus_{\alpha \in \Delta_\Sigma} \mfrak{g}_\alpha,
\end{align}
and the nilradical $\mfrak{u}$ of $\mfrak{p}$ and the opposite nilradical $\widebar{\mfrak{u}}$ are
\begin{align}
  \mfrak{u}= \bigoplus_{\alpha \in \Delta^+ \smallsetminus \Delta_\Sigma^+}\mfrak{g}_\alpha, \qquad
	\widebar{\mfrak{u}}=
	\bigoplus_{\alpha \in \Delta^+ \smallsetminus \Delta_\Sigma^+} \mfrak{g}_{-\alpha},
\end{align}
respectively. The $\Sigma$-height $\htt_\Sigma(\alpha)$ of a root $\alpha \in \Delta$ is defined by
\begin{align}
  \htt_\Sigma(a_1\alpha_1+a_2\alpha_2) = a_1,
\end{align}
and $\mfrak{g}$ is a $|1|$-graded Lie algebra with respect to the grading given by $\mfrak{g}_i = \bigoplus_{\alpha \in \Delta,\, \htt_\Sigma(\alpha)=i} \mfrak{g}_\alpha$ for $0 \neq i \in \Z$, and $\mfrak{g}_0= \mfrak{h} \oplus \smash{\bigoplus}_{\alpha \in \Delta,\, \htt_\Sigma(\alpha)=0}
\mfrak{g}_\alpha$. In particular, we have $\mfrak{u}=\mfrak{g}_1\simeq \C^2$,
$\widebar{\mfrak{u}}=\mfrak{g}_{-1}\simeq\C^2$ and
$\mfrak{l}= \mfrak{g}_0 \simeq \C \oplus \mfrak{sl}(2,\C)$.

The basis $\{e_1, e_2\}$ of the root spaces in the nilradical
$\mfrak{u}$ is given by
\begin{align}
  e_1 = \begin{pmatrix}
    0 & 1 & 0  \\
    0 & 0 & 0  \\
    0 & 0 & 0
  \end{pmatrix}\!, \qquad
  e_2 = \begin{pmatrix}
    0 & 0 & 1 \\
    0 & 0 & 0 \\
    0 & 0 & 0
  \end{pmatrix}\!,
\end{align}
the basis $\{f_1, f_2\}$ of the root spaces in the opposite nilradical
$\widebar{\mfrak{u}}$ is
\begin{align}
  f_1 = \begin{pmatrix}
    0 & 0 & 0  \\
    1 & 0 & 0  \\
    0 & 0 & 0
  \end{pmatrix}\!, \qquad
  f_2 = \begin{pmatrix}
    0 & 0 & 0 \\
    0 & 0 & 0 \\
    1 & 0 & 0
  \end{pmatrix}\!,
\end{align}
and finally the basis $\{h_0, h, e, f\}$ of the Levi subalgebra $\mathfrak{l}$ is
\begin{align}
\begin{gathered}
 e  = \begin{pmatrix}
    0 & 0 & 0  \\
    0 & 0 & 1  \\
    0 & 0 & 0
  \end{pmatrix}\!, \quad
   h = \begin{pmatrix}
    0 & 0 & 0 \\
    0 & 1 & 0 \\
    0 & 0 & -1
  \end{pmatrix}\!, \quad
  f = \begin{pmatrix}
    0 & 0 & 0 \\
    0 & 0 & 0 \\
    0 & 1 & 0
  \end{pmatrix}\!, \\
  h_0 = \begin{pmatrix}
    1 & 0 & 0  \\
    0 & -\frac{1}{2} & 0  \\
    0 & 0 & -\frac{1}{2}
  \end{pmatrix}\!,
\end{gathered}
\end{align}
where $h_0$ generates a basis of the center $\mfrak{z}(\mfrak{l})$ of $\mfrak{l}$.

Any character $\chi \in \Hom_P(\mfrak{p},\C)$ is given by
\begin{align}
\chi= \lambda \widetilde{\omega}_1, \quad \lambda \in \C ,
\end{align}
where $\widetilde{\omega}_1\in \Hom_P(\mfrak{p},\C)$
is equal to $\omega_1\in \mfrak{h}^*$ regarded as trivially extended to $\mfrak{p}=\mfrak{h} \oplus \mfrak{g}_{\alpha_2} \oplus \mfrak{g}_{-\alpha_2} \oplus \mfrak{u}$. Throughout the article we use the 
simplified notation $\lambda$ for $\lambda \widetilde{\omega}_1$.
The vector $\rho \in \Hom_P(\mfrak{p},\C)$ defined by the formula \eqref{eq:rho vector} is then
\begin{align}
\rho={\textstyle {3 \over 2}}\,\widetilde{\omega}_1.
\end{align}


\section{$\widetilde{\mathrm{SL}}(3,\mR)$ and the symplectic Dirac operator}
\label{sec:sl3sdo}

In this section we retain the notation in Section \ref{sec:projstr}
and describe the class of representations of $\mathfrak{g}$
on the space of sections of vector bundles on
$\smash{\widetilde{G}}/\smash{\widetilde{P}}$ associated to the
simple metaplectic submodules of the Segal-Shale-Weil representation
$\mathbb{S}_{\lambda+\rho}$ of $\smash{\widetilde{P}}$ twisted by
characters $\lambda + \rho \in \Hom_{\widetilde{P}}(\mfrak{p},\C)$.

The induced representations in question are described in
the non-compact picture, given by restricting sections to
the open Schubert cell $U_e\subset \widetilde{G}/\widetilde{P}$ which is 
isomorphic by the
exponential map to the opposite nilradical $\widebar{\mfrak{u}}(\mR)$.
We denote by $(\hat{x}, \hat{y})$
the linear coordinate functions on
$\widebar{\mfrak{u}}(\mR)$ with respect to the basis $\{f_1, f_2\}$ of $\widebar{\mfrak{u}}(\mR)$,
and by $(x, y)$ the dual linear coordinate functions on
$\widebar{\mfrak{u}}^*(\mR)$. The Weyl algebra
$\eus{A}^\mfrak{g}_{\widebar{\mfrak{u}}}$ is generated by
\begin{align}
\{\hat{x}, \hat{y}, \partial_{\hat{x}}, \partial_{\hat{y}}\}
\label{eq:Weyl algebra generators}
\end{align}
and the Weyl algebra
$\eus{A}^\mfrak{g}_{\widebar{\mfrak{u}}^*}$ by
\begin{align}
\{x, y, \partial_{x}, \partial_{y}\}. \label{eq:dual Weyl algebra generators}
\end{align}

For a $\mfrak{p}$-module $(\sigma,\mathbb{V})$, $\sigma \colon \mfrak{p} \rarr \mfrak{gl}(\mathbb{V})$,
the twisted $\mfrak{p}$-module $(\sigma_\lambda, \mathbb{V}_\lambda)$,
$\sigma_\lambda \colon \mfrak{p} \rarr \mfrak{gl}(\mathbb{V}_\lambda)$, with a twist
$\lambda \in \Hom_{\widetilde{P}}(\mfrak{p},\C)$ is defined by
\begin{align}
  \sigma_\lambda(X)v=\sigma(X)v+\lambda(X)v
\end{align}
for all $X \in \mfrak{p}$ and $v \in \mathbb{V}_\lambda \simeq \mathbb{V}$ (the isomorphism of vector spaces).
\medskip

We use the following realization of the simple $\mathfrak{mp}(2,\C)$-submodules of the
Segal-Shale-Weil representation. The Fock model is the unitarizable $\mfrak{mp}(2,\C)$-module
$\mathbb{S}=\C[q]$. Since the simple part of the Levi algebra
$\mathfrak{l}$ is $\mfrak{l}^{\rm s} \simeq \mfrak{mp}(2,\C)$, we
realize the simple metaplectic submodules of the Segal-Shale-Weil representation as
the representations of $\mfrak{l}^{\rm s}$ on the subspace of polynomials of even and
odd degree, respectively. The generators act as
\begin{align}\label{eq:segal-shale-weil}
   \sigma(e)= {\textstyle {\imag \over 2}} \partial^2_{q}, \qquad \sigma(h)=-q\partial_{q}-{\textstyle {1 \over 2}}, \qquad
 \sigma(f)= {\textstyle {\imag \over 2}} q^2.
\end{align}
The scalar product $\langle \cdot\,,\cdot \rangle \colon \mathbb{S} \otimes_\C \mathbb{S} \rarr \C$ on $\mathbb{S}$ is defined through the $\mfrak{l}^{\rm s}$-equivariant embedding into the space of Schwartz functions
$\iota \colon \mathbb{S} \rarr \mcal{S}(\R)$,
\begin{align}
  \langle p_1,p_2 \rangle = \int_{\mR} \iota(\widebar{p}_1)\iota(p_2)\,{\rm d}q
	\quad \text{for all}\quad p_1,p_2\in \mathbb{S}.
\end{align}
The representation of $\mfrak{l}^{\rm s}$ is then extended to a representation
of $\mathfrak{p}$ by the trivial action of the center $\mfrak{z}(\mfrak{l})$
of $\mfrak{l}$ and by the trivial action of the nilradical $\mathfrak{u}$ of
$\mfrak{p}$. We retain the same notation $\sigma \colon \mfrak{p} \rarr \mfrak{gl}(\mathbb{S})$
for the extended action of the parabolic subalgebra
$\mfrak{p}$ of $\mfrak{g}$. In what follows, we are interested in the twisted $\mfrak{p}$-module
$\sigma_\lambda \colon \mfrak{p} \rarr \mfrak{gl}(\mathbb{S}_\lambda)$ with a twist
$\lambda \in \Hom_{\widetilde{P}}(\mfrak{p},\C)$.
\medskip

\theorem{\label{reproper}Let $\lambda \in \Hom_{\widetilde{P}}(\mfrak{p},\C)$. Then the embedding of $\mfrak{g}$ into
$\eus{A}^\mfrak{g}_{\widebar{\mfrak{u}}} \otimes_\C \End \mathbb{S}_{\lambda+\rho}$ and
$\eus{A}^\mfrak{g}_{\widebar{\mfrak{u}}^*}\! \otimes_\C \End \mathbb{S}_{\lambda+\rho}$ is given by
\begin{enumerate}
\item[1)]
\begin{align}
  \pi_\lambda(f_1)=-\partial_{\hat{x}}, \qquad \pi_\lambda(f_2)=-\partial_{\hat{y}},
\end{align}
\begin{align}
  \hat{\pi}_\lambda(f_1)=-x, \qquad \hat{\pi}_\lambda(f_2)=-y;
\end{align}
\item[2)]
\begin{align}
\begin{gathered}
 \pi_\lambda(e)=-\hat{y} \partial_{\hat{x}}+ {\textstyle \frac{\imag}{2}} \partial_q^2,\quad
 \pi_\lambda(h)=- \hat{x} \partial_{\hat{x}}+\hat{y}\partial_{\hat{y}}-q \partial_q-{\textstyle \frac{1}{2}}, \quad
 \pi_\lambda(f)=- \hat{x} \partial_{\hat{y}}+ {\textstyle \frac{\imag}{2}} q^2, \\
 \pi_\lambda(h_0)={\textstyle {3 \over 2}}(\hat{x} \partial_{\hat{x}}+\hat{y}\partial_{\hat{y}})+\lambda +{\textstyle \frac{3}{2}},
\end{gathered}
\end{align}
\begin{align}
\begin{gathered}
 \hat{\pi}_\lambda(e) = x\partial_y+{\textstyle \frac{\imag}{2}}\partial_q^2, \quad \hat{\pi}_\lambda(h) = x\partial_x-y\partial_y-q\partial_q-{\textstyle \frac{1}{2}}, \quad
\hat{\pi}_\lambda(f) = y\partial_x+{\textstyle \frac{\imag}{2}} q^2, \\
\hat{\pi}_\lambda(h_0) = -{\textstyle \frac{3}{2}}(x\partial_x+y\partial_y)+\lambda-{\textstyle \frac{3}{2}};
\end{gathered}
\end{align}
\item[3)]
\begin{align}
\begin{aligned}
{\pi}_\lambda(e_1) &=\hat{x}(\hat{x}\partial_{\hat{x}}+\hat{y}\partial_{\hat{y}}+\lambda+{\textstyle {3\over 2}}) + {\textstyle {1\over 2}}\hat{x}\big(q\partial_q+{\textstyle {1 \over 2}}\big) -{\textstyle {\imag \over 2}}\hat{y}q^2, \\
 {\pi}_\lambda(e_2) &=\hat{y}(\hat{x}\partial_{\hat{x}}+\hat{y}\partial_{\hat{y}}+\lambda+{\textstyle {3\over 2}})- {\textstyle {1\over 2}}\hat{y}\big(q\partial_q+{\textstyle {1 \over 2}}\big) -{\textstyle {\imag \over 2}}\hat{x}\partial_q^2,
\end{aligned}
\end{align}
\begin{align}
\begin{aligned}
\hat{\pi}_\lambda(e_1) &=\partial_x\big(x\partial_x+y\partial_y-\lambda+{\textstyle \frac{1}{4}}\big) +{\textstyle \frac{1}{2}}q(\imag q\partial_y-\partial_x\partial_q), \\
\hat{\pi}_\lambda(e_2) &=\partial_y\big(x\partial_x+y\partial_y-\lambda+{\textstyle \frac{1}{4}}\big) -{\textstyle \frac{\imag}{2}} \partial_q(\imag q\partial_y-\partial_x\partial_q).
\end{aligned}
\end{align}
\end{enumerate}}

\proof{The proof of this claim for the trivial representation of $\mfrak{p}$ instead
of the simple metaplectic submodules of the Segal-Shale-Weil representation follows from 
Theorem $1.3$ in \cite{Krizka-Somberg2015}. For the Segal-Shale-Weil representation it follows by a straightforward verification of all commutation relations for the Lie algebra $\mfrak{g}$.}

\theorem{\label{xssingvect}The space $\Sol(\mfrak{g},\mfrak{p};\C[\widebar{\mfrak{u}}^*]
\otimes_\C \mathbb{S}_{\lambda+\rho})^\mcal{F}$ is for
$\lambda=\frac{3}{4}\widetilde{\omega}_1$ non-trivial and contains the
$L$-submodule $X_s\mathbb{S}_{\lambda+\rho}$.}

\proof{By Theorem \eqref{reproper}, the $\mathfrak{sl}(2,\C)$-algebra
structure for the operators \eqref{sl2howedual} implies
\begin{align*}
\hat{\pi}_\lambda(e_1)X_sv_0 &=
\partial_x\big(1-\lambda+{\textstyle \frac{1}{4}}\big)X_sv_0+
{\textstyle \frac{1}{2}}q[D_s,X_s]v_0  = \big({\textstyle \frac{5}{4}}-\lambda\big)[\partial_x,X_s]v_0
-{\textstyle \frac{\imag}{2}}qv_0 \\
& = \big({\textstyle \frac{5}{4}-\lambda-\frac{1}{2}}\big)\imag qv_0 =0,
\\
\hat{\pi}_\lambda(e_2)X_sv_0 &=
\partial_y\big(1-\lambda+{\textstyle \frac{1}{4}}\big)X_sv_0
-{\textstyle \frac{1}{2}}\imag \partial_q[D_s,X_s]v_0 = \big({\textstyle \frac{5}{4}}-\lambda\big)[\partial_y,X_s]v_0
-{\textstyle \frac{1}{2}}\partial_q v_0 \\
& = \big({\textstyle \frac{5}{4}-\lambda-\frac{1}{2}}\big)\partial_qv_0 =0
\end{align*}
for all $v_0 \in \mathbb{S}_{\lambda+\rho}$, provided $\lambda=\frac{3}{4}$.
The proof is complete.}

\remark{We notice that the vectors $X_s^kv$ for $k > 1$
are not singular vectors. For example, a straightforward computation
reveals
\begin{align}
\hat{\pi}_\lambda(e_1)X_s^kv_0 &=
\big(k-\lambda+{\textstyle \frac{1}{4}}\big)\big(\imag kqX_s^{k-1}
+\imag {\textstyle \frac{k(k-1)}{2}}yX_s^{k-2}\big)v
-\imag {\textstyle \frac{k(k-1)}{4}}qX_s^{k-1}v,
\end{align}
which is non-zero for all $\lambda\in\mC$.}

It is not difficult to exploit the results in \cite{DeBie-Somberg-Soucek2014}
in order to classify the complete set of solutions of the system in
Theorem \ref{reproper}, but the detailed analysis goes beyond the
scope of our article.

Theorem \ref{xssingvect} has the following classical corollary, which
explains the relationship between the geometrical problem of finding
$G$-equivariant differential operators between induced representations
and the algebraic problem of finding homomorphisms between generalized
Verma modules, cf.\ \cite{Collingwood-Shelton1990}, \cite{Krizka-Somberg2016}.
There is a double cover $\smash{\widetilde{P}}=(\GL(1,\mR)_+\times\smash{\widetilde{\SL}}(2,\mR))\ltimes \mR^2$
of the maximal parabolic subgroup
$P$ of the Lie group $\SL(3,\mR)$, which splits over the unipotent subgroup $N\simeq \mR^2$
in the Langlands-Iwasawa decomposition of $\smash{\widetilde{P}}$, \cite{Wolf1976}. Let us
note that the extension cocycle splits over the field of complex numbers.
\medskip

\theorem{\label{thm:edo}Let $\smash{\widetilde{G}}=\smash{\widetilde{\SL}}(3,\mR)$ and let $\smash{\widetilde{P}}=(\GL(1,\mR)_+\!\times\smash{\widetilde{\SL}}(2,\mR))\ltimes \mR^2$
be the maximal parabolic subgroup of $\smash{\widetilde{G}}$, whose unipotent
subgroup in the Langlands-Iwasawa decomposition of $\smash{\widetilde{P}}$
is $N\simeq \mR^2$. For $\mathbb{V}=\mathbb{S}_\lambda$ we have $\mathbb{V}^* \simeq \mathbb{S}^*_{-\lambda}$.
Then the singular vector constructed in Theorem \ref{xssingvect} corresponds
to the $\smash{\widetilde{G}}$-equivariant differential operator, given
in the non-compact picture of the induced representations by
\begin{align}
\begin{gathered}\label{projformds}
D_s \colon {\lC}^\infty(\widebar{\mfrak{u}}(\R),\mS^*_{{3\over 4}\widetilde{\omega}_1})
\rarr  {\lC}^\infty(\widebar{\mfrak{u}}(\R),\mS^*_{{9\over 4}\widetilde{\omega}_1}),\\
\varphi \mapsto (\imag q\partial_{\hat{y}}-\partial_{\hat{x}}\partial_{q})\varphi.
\end{gathered}
\end{align}
The infinitesimal intertwining property of $D_s$ is
\begin{align}
  D_s \pi^*_{-{3 \over 4}}(X)=\pi^*_{{3 \over 4}}(X)D_s
\end{align}
for all $X \in \mfrak{sl}(3,\R)$.
With abuse of notation we used the same symbol $D_s$ and the
same terminology "the symplectic Dirac operator" as in \eqref{sl2howedual}
due to the coincidence of \eqref{sl2howedual} and \eqref{projformds}.}

Let us finally explain the notion of the dual representation $\mS^*$.
Let us define a non-degenerate pairing
$(\cdot\,,\cdot ) \colon \mathbb{S} \otimes_\C \mathbb{S} \rarr \C$ on $\mathbb{S}$
by the formula
\begin{align}
  (p_1(q),p_2(q)) = p_1(\partial_q)p_2(q)|_{q=0}.
\end{align}
Then we can identify the (restricted) dual space to $\mathbb{S}$ with $\mathbb{S}$ and
the structure of the dual $\mfrak{p}$-module on $\mathbb{S}$ is given as follows:
the generators of $\mfrak{mp}(2,\C)$ act on $\mS^*$ by
\begin{align}
\begin{aligned} \label{eq:segal-shale-weil dual}
  \sigma^*(e)= -{\textstyle {\imag \over 2}} q^2, \qquad \sigma^*(h)=q\partial_{q}+{\textstyle {1 \over 2}}, \qquad
 \sigma^*(f)= -{\textstyle {\imag \over 2}}\partial_{q}^2,
\end{aligned}
\end{align}
while the generator $h_0$ of the center $\mfrak{z}(\mfrak{l})\subset\mfrak{l}$
and of the nilradical $\mathfrak{u}$ of $\mfrak{p}$ act trivially. We notice 
that this representation is compatible with \eqref{HXYxy}.


\section*{Acknowledgement}

The authors gratefully acknowledge the support of the grant GA\,CR P201/12/G028 and SVV-2016-260336.



\providecommand{\bysame}{\leavevmode\hbox to3em{\hrulefill}\thinspace}
\providecommand{\MR}{\relax\ifhmode\unskip\space\fi MR }
\providecommand{\MRhref}[2]{%
  \href{http://www.ams.org/mathscinet-getitem?mr=#1}{#2}
}
\providecommand{\href}[2]{#2}

\end{document}